\documentclass[onefignum,onetabnum]{siamart190516}

\usepackage{lipsum}
\usepackage{amsfonts}
\usepackage{graphicx}
\usepackage{epstopdf}
\usepackage{algorithmic}
\ifpdf
  \DeclareGraphicsExtensions{.eps,.pdf,.png,.jpg}
\else
  \DeclareGraphicsExtensions{.eps}
\fi

\newsiamremark{remark}{Remark}
\newsiamremark{hypothesis}{Hypothesis}
\crefname{hypothesis}{Hypothesis}{Hypotheses}
\newsiamthm{claim}{Claim}

\headers{Reparametrization of Curves via Interpolation}{Kazuki Koga}

\title{Numerical Reparametrization of Periodic Planar Curves via Curvature Interpolation\thanks{
\funding{This work was partially funded by Kyoto University Grant for Start-Up 2021, JST ACT-X JPMJAX2106, and JSPS KAKENHI Grant No.\,21K20325. }}}

\author{Kazuki Koga\thanks{\email{koga.kazuki.2m@kyoto-u.ac.jp}.}}

\usepackage{amsopn}

\ifpdf
\hypersetup{
  pdftitle={Numerical Reparametrization of Periodic Planar Curves via Curvature Interpolation},
  pdfauthor={K. Koga}
}
\fi

\begin{document}

\maketitle

\begin{abstract}
A novel static algorithm is proposed for numerical reparametrization of periodic planar curves. The method identifies a monitor function of the arclength variable with the true curvature of an open planar curve and considers a simple interpolation between the object and the unit circle at the curvature level. Since a convenient formula is known for tangential velocity that maintains the equidistribution rule with curvature-type monitor functions, the strategy enables to compute the correspondence between the arclength and another spatial variable by evolving the interpolated curve. With a certain normalization, velocity information in the motion is obtained with spectral accuracy while the resulting parametrization remains unchanged. Then, the algorithm extracts a refined representation of the input curve by sampling its arclength parametrization whose Fourier coefficients are directly accessed through a simple change of variables. As a validation, improvements to spatial resolution are evaluated by approximating the invariant coefficients from downsampled data and observing faster global convergence to the original shape.
\end{abstract}



\section{Introduction}\label{intro}
Parametric planar curves are a basic apparatus for representing the moving boundaries of two-dimensional and axisymmetric objects. Typical settings include the Lagrangian description of fluid interfaces and their numerical simulations, where control of parametrization has major advantages beyond stable spatial discretization. For instance, the uniform parametrization is used for simplifying the source of stiffness due to surface tension\,\cite{HoLoSh1994}, for reducing the size of unknowns in searching for time-periodic interface motion based on numerical optimization\,\cite{AmWi2010}, and for computing reference solutions in a geometric convergence study of simulations of axisymmetric droplets\,\cite{Koga2021}. On the other hand, in numerical studies of singularity formations such as capillary pinch-off\,\cite{LeLi2003}, locality of parametrization is key to resolving subtle behavior at small length-scales. \par
A common approach to enforcing a specific parametrization is to exploit the arbitrariness of tangential velocity. That is, noting that the shape of a moving boundary is completely determined by normal velocity, as rigorously proved for the mean curvature flow in the plane\,\cite{EpGa1987}, the tangential counterpart can be defined to adjust its parametrization and hence a mesh distribution after spatial discretization. Among those, the uniform parametrization is a typical choice for various moving boundary problems, not only because of the reasons above but because the corresponding tangential velocity is computed by simple numerical anti-derivatives\,\cite{HoLoSh1994}. An analogous situation is found by the same authors for the equidistribution rule with curvature-type criterions, which also leads to solvable equations for tangential velocity and has been applied to several non-singular problems\,\cite{SeYa2011, Kropinski2001}. On the other hand, Nitsche and Steen\,\cite{NiSt2004} suggests a dynamical strategy for constructing parametrizations with non-geometric guidelines to resolve the singularity of axisymmetric capillary pinch-off. Their simulations start from the uniform mesh and advect discrete points continuously toward a prescribed distribution in a simpler framework than the equidistribution rule.
The core idea in the previous work is further extended by Koga\,\cite{Koga2021}, which introduces signal processing for designing effective guidelines with far less human intervention. In those schemes, however, the use of non-trivial criterions leads to practically unsolvable equations for tangential velocity as part of solutions, and requires a backward temporal approximation \cite{NiSt2004} that can cause some order reduction of time integrators, as experimentally demonstrated in \cite{Koga2021}. Moreover, resulting tangential motion is smooth in time only if signal processing is built upon extremely well-behaved operators.
\par
In practical simulations, a static reparametrization algorithm allows to perform mesh refinement and time-stepping independently. Successful applications are, for example, the work by Leppinen and Lister\,\cite{LeLi2003} on axisymmetric capillary pinch-off, which advects computational points by the Lagrangian velocity \cite{CaLi1992} and redistributes them at each time step so that the mesh density remains proportional to the distance from the shrinking necks. A similar strategy is used by Burton and Taborek \cite{BuTa2007} to compare breakup of two-dimensional droplets against the axisymmetric case. Although this static approach does not have a solid mathematical foundation and only works under strong assumptions on initial conditions, it is expected to overcome the difficulties of the dynamical method above and further extend the recent work\,\cite{Koga2021}. In this direction, Seol and Lai \cite{SeLa2020} proposes an optimization-based algorithm for statically achieving the equidistribution rule \cite{HoLoSh1994} with smooth complexity measures. It can be regarded as a spectrally accurate version of {de Boor's} algorithm \cite{Boor1973} whose convergence is theoretically analyzed by Xu et al.\,\cite{XuHuRuWi2011}, but requires proper initial guesses to converge and is sensitive to oscillatory perturbations on input data. \par
As a robust alternative to the iterative algorithm\,\cite{SeLa2020}, the present work develops a simulation-based approach to the static reparametrization problem. In our method, a monitor function of the arclength variable, which describes local complexity of a planar curve by strictly positive values, is regarded as the true curvature of another curve in the canonical parametrization and linearly interpolated with that of the unit circle as the other end-point. This one-parameter family of functions defines nonlinear motion of a single planar curve without forming a cusp, and its time evolution is completely known from the interpolation at the curvature level. Then, due to the simplicity of tangential velocity that maintains the equidistribution rule with the positive curvature\,\cite{HoLoSh1994}, the desired parametrization (or equivalently the mapping of the computational domain onto the arclength interval) is obtained by simulating the artificial time evolution by a time-stepping method. Moreover, a $L^1$-type normalization of monitor functions, which periodizes derivatives of the arclength parametrization of the moving curve, guarantees spectrally accurate discretization in space while the resulting parametrizations remain exactly the same. Once a new parametrization is found, the refined representation of the input curve is computed via its arclength paramerization as intermediate data whose Fourier coefficients are directly accessed from the initial parametrization\,\cite{Koga2021}. As a result, the algorithm involves the forward and inverse Fourier transforms with unevenly distributed nodes, and those operations are accelerated by a modern implementation of the nonuniform fast Fourier transform (NUFFT) \cite{BaMaKl2019}. In numerical experiments, where we consider two nearly singular examples with fluid interfaces in mind, improvements to spatial resolution are confirmed by recomputing the Fourier coefficients of the arclength parametrization from downsampled data and observing faster global convergence to the original shapes. To our best knowledge, such a quantitative aspect has been overlooked in the literature due to the lack of invariants that enable to reconstruct a planar curve in a canonical form, and therefore constitutes the main contributions of the present work.


The rest of this paper is organized as follows. Section \ref{formulate} formulates the equidistribution rule for planar curves and the one-parameter family of curves that arises from the curvature interpolation. Section \ref{method} explains our numerical methods and some technical issues in their efficient implementations. Section \ref{results} presents numerical results and validates the implemented codes. Concluding remarks are given in Section \ref{conclusion}.


\section{Formulations}\label{formulate}
This section formulates the equidistribution rule for planar curves and the one-parameter family of curves that arises from linearly interpolating two periodic curvatures. In the following, subscripts such as $\alpha$, $s$ and $t$ denote differentiation with respect to the corresponding variables, and the operators $\partial^{-1}_\alpha$ and $\partial^{-1}_s$ are the anti-derivatives with the base points $\alpha=0$ and $s=0$, respectively.
\subsection{Elementary differential geometry}
Let $\mathbf{X}: [0,2\pi]\times [0,1] \rightarrow \mathbb{R}^2$ be a one-parameter family of smooth planar curves parametrized by the variable $\alpha$:
\begin{equation}
\label{eq:def_prmt}
\mathbf{X}(\alpha,t) = (x(\alpha,t), y(\alpha,t)),\quad   \alpha \in [0,2\pi],\quad  t \in [0,1].
\end{equation}
That is, for each $t$, the mapping $\mathbf{X}$ represents a planar curve whose coordinates $x$ and $y$ are both smooth in $\alpha$ and satisfy $\sqrt{x_\alpha^2 + y_\alpha^2}>0$. In particular, a smooth planar curve is called periodic if its parametrization can be extended as a $2\pi$-periodic mapping. We refer to the function $s_\alpha=\sqrt{x_\alpha^2 + y_\alpha^2}$ as the {\it  local spacing}, and it is one of the central concepts in the reparametrization problem. \par
We identify a parametric planar curve with another if their images are the same set in $\mathbb{R}^2$. Of such representations, most important is the arclength parametrization that associates points on the curve with the arclength $s$ measured from a base point. Assuming that the mapping (\ref{eq:def_prmt}) is bijective, the arclength $s$ is defined to be a monotonically increasing function of $\alpha$ via integrating its derivative $s_\alpha$:
\begin{equation}
\label{eq:def_arclength}
s(\alpha,t)=[\partial^{-1}_\alpha s_\alpha](\alpha,t).
\end{equation}
Therefore, the arclength parametrization $\tilde{\mathbf{X}}$ is written as
\begin{equation}
\label{eq:X_arclength}
\tilde{\mathbf{X}}(s,t) =\mathbf{X}(\alpha(s,t),t),\quad   s \in [0,L],
\end{equation}
where $L=s(2\pi,t)$ is the total length of the curve and the function $\alpha$ is the inverse of (\ref{eq:def_arclength}). Conversely, the original mapping (\ref{eq:def_prmt}) is regarded as the composition
\begin{equation}
\label{eq:X_composite}
\mathbf{X}(\alpha,t) =\tilde{\mathbf{X}}(s(\alpha,t),t),\quad   \alpha \in [0,2\pi],
\end{equation}
and thus a change of paramerizations is reduced to finding a new one-to-one correspondence between the two intervals $[0,2\pi]$ and $[0,L]$. We accomplish this task by constructing the local spacing $s_\alpha$ rather than its primitive function $s$.  \par
\begin{figure}[h]
\begin{center}
 \includegraphics[width=0.4\linewidth]{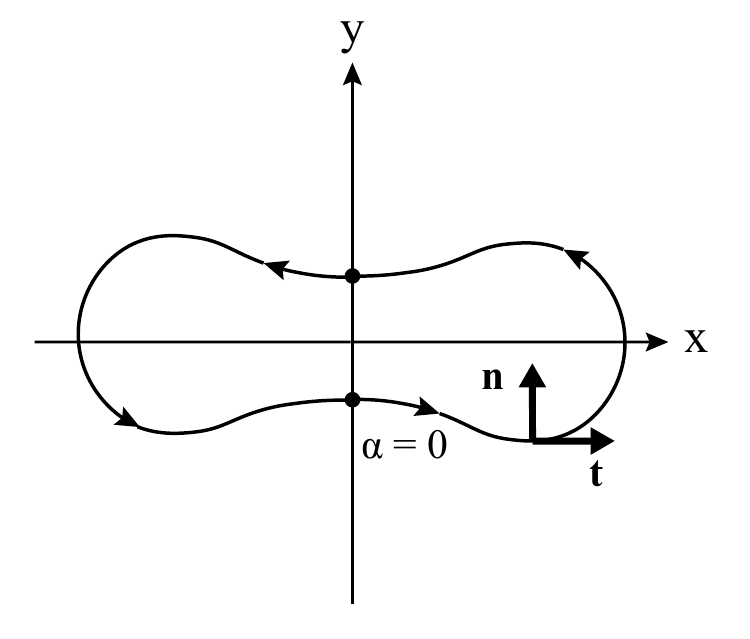}
\end{center}
\caption{\label{fig:fig_plane_curve} Parametrization of planar curve and corresponding pair of normal vectors $(\mathbf{n},\mathbf{t})$. } 
\end{figure}

Typically, the second variable $t$ is interpreted as time and the evolution equation for the parametric curve $\mathbf{X}$ is derived in the form 
\begin{equation}
\label{eq:decompX}
 \mathbf{X}_t = U\mathbf{n}+V\mathbf{t}.
\end{equation}
Here, $\mathbf{n}$ and $\mathbf{t}$ are the unit normal and tangent vectors satisfying the Frenet-Serret formula $\mathbf{t}_s=\kappa\mathbf{n}$ with the signed curvature
\begin{equation}
\label{eq:curvature}
\kappa={y_{ss}x_{s}-y_{s}x_{ss}},
\end{equation}
where the operator $\frac{\partial}{\partial s}$ is defined as $\frac{\partial}{\partial s}=s_\alpha^{-1} \frac{\partial}{\partial \alpha}$. More specifically, we choose
\begin{equation}
\label{eq:def_unit_nt}
\mathbf{n}=(-y_s,x_s),\quad \mathbf{t}=(x_s,y_s),
\end{equation}
so that the unit circle has the curvature $\kappa \equiv 1$ with a counterclockwise parametrization (see Fig.\,\ref{fig:fig_plane_curve}). 
The orthogonal decomposition (\ref{eq:decompX}) plays an important role in controlling the parametrization, because the shape of a moving planar curve is determined solely by the normal velocity $U$, and the tangential velocity $V$ only affects the local spacing $s_\alpha$. Following {Hou et al.\,\cite{HoLoSh1994}}, we switch from the Cartesian coordinates $(x,y)$ to the so-called angle--arclength variables $(\theta,s_\alpha)$ defined via 
\begin{equation}
\label{eq:def_theta}
\mathbf{X}_\alpha =(s_\alpha\cos\theta, s_\alpha \sin \theta).
\end{equation}
Therefore, the variable $\theta$ measures the angle between the $x$-axis and the vector $\mathbf{t}$, and it is directly linked to the curvature $\kappa$ as $\theta_s =\kappa$. With this property, the evolution equations for the new variables $(\theta,s_\alpha)$ follow from the decomposition (\ref{eq:decompX}) and the Frenet-Serret formula $\mathbf{t}_s=\kappa\mathbf{n}$: 
 \begin{equation}
 \label{eq:evol_theta_sa}
    \theta_t = \frac{U_\alpha + V \theta_\alpha}{s_\alpha}, \quad s_{\alpha,t} = V_\alpha - \theta_\alpha U,
  \end{equation}
As mentioned earlier, because of the tangential velocity being arbitrary, it is possible to control the local spacing $s_\alpha$ while keeping the image of $\mathbf{X}$ unchanged. Moreover, since we are able to reconstruct $(x,y)$ up to a translation by integrating Eq.\,(\ref{eq:def_theta}), the angle-arclength representation is considered most natural in defining motion of a planar curve by interpolating two curvatures and extracting the local spacing $s_\alpha$ from the moving object.
\subsection{Equidistribution rule}
Let $\varphi$ be a strictly positive function of the arclength variable $s$ and define the local spacing $s_\alpha$ by the equation
      \begin{equation}
 \label{eq:def_equidist}
      s_\alpha \varphi = \frac{1}{2\pi} \int_0^{2\pi} s_\alpha \varphi d\alpha,
     \end{equation}
where $\varphi$ above is understood as $\varphi(\alpha) = \varphi(s(\alpha))$. The parametrization (\ref{eq:def_equidist}) is called the {\it equidisitribution rule} for planar curves, and the user-defined factor $\varphi$ is called the {\it monitor function}, which describes complexity of an object in the arclength parametrization. The equidistribution rule is originally introduced in the context of boundary-value problems of ordinary differential equations \cite{Boor1973}, and similar concepts have been sought for higher dimensions in a variational form \cite{ReWa2000}. The variant (\ref{eq:def_equidist}) is first suggested by Hou et al.\,\cite[Appendix 2]{HoLoSh1994} as an alternative to the uniform parametrization (i.e. $\varphi \equiv 1$), and they show that curvature-dependent monitor functions $\varphi(\kappa)$ (e.g., $\varphi = 1+\kappa^2$) lead to a convenient formula for tangential velocity. To see this, from Eqs.\,(\ref{eq:evol_theta_sa}), differentiating both sides of $\kappa = (\theta_\alpha/s_\alpha)$ with respect to $t$ yields
  \begin{equation}
    \kappa_t = \frac{1}{s_\alpha}\biggl\{\biggl(\frac{U_\alpha + V \theta_\alpha}{s_\alpha}\biggr)_\alpha-\kappa(V_\alpha - \theta_\alpha U) \biggr\},
  \end{equation}
and hence, with the notation $\varphi_\kappa=\frac{d \varphi}{d \kappa}$, we have
\begin{align}
\nonumber (s_\alpha \varphi)_t &= s_{\alpha,t}\varphi +s_\alpha \kappa_t \varphi_\kappa \\
\nonumber                                 &=(V_\alpha - \theta_\alpha U)\varphi+ \biggl\{\biggl(\frac{U_\alpha + V \theta_\alpha}{s_\alpha}\biggr)_\alpha-\kappa(V_\alpha - \theta_\alpha U) \biggr\}\varphi_\kappa\\
\label{eq:def_diff_sa_phi}          &= (\varphi V)_\alpha+\varphi_\kappa (U_s)_\alpha+\theta_\alpha U(\kappa \varphi_\kappa-\varphi).
\end{align}
Then, plugging (\ref{eq:def_diff_sa_phi}) into the definition  (\ref{eq:def_equidist}), a linear integro-differential equation for the component $V$ is obtained as
\begin{align}
\label{eq:diff_Vphi}
 (\varphi V)_\alpha&= -[\varphi_\kappa (U_s)_\alpha+\theta_\alpha U(\kappa \varphi_\kappa-\varphi)]\\
\nonumber            &+ \frac{1}{2\pi} \int_0^{2\pi} [(\varphi V)_\alpha+\varphi_\kappa (U_s)_\alpha+\theta_\alpha U(\kappa \varphi_\kappa-\varphi)] d\alpha.
 \end{align}
 To solve Eq.\,(\ref{eq:diff_Vphi}), we need boundary values for $V$ that are consistent with physical velocity at the end-points of a moving curve.  For example, the boundary values
 \begin{equation}
\label{eq:def_Vend_period}
V(0,t)=V(2\pi,t)=0,
\end{equation}
often arise in axisymmetric problems \cite{Koga2021} and yield an explicit integral form 
 \begin{align}
\label{eq:int_V} V(\alpha,t)&= -\frac{1}{\varphi}\partial^{-1}_\alpha [\varphi_\kappa (U_s)_\alpha+\theta_\alpha U(\kappa \varphi_\kappa-\varphi)]\\
\nonumber            &+\frac{\alpha}{2\pi\varphi} \int_0^{2\pi} [\varphi_\kappa (U_s)_\alpha+\theta_\alpha U(\kappa \varphi_\kappa-\varphi)] d\alpha.
 \end{align}
 Later, we derive a special case of (\ref{eq:diff_Vphi}) with the monitor function $\varphi \equiv \kappa$, which does not involve integrals and can be solved by algebraic operations on the known data.
\subsection{Curvature interpolation}
As pointed out in \cite{HoLoSh1997}, any local spacing $s_\alpha$ can be written as a product of the total length $L$ and a strictly positive function $R$:
\begin{equation}
\label{eq:def_R}
s_\alpha(\alpha,t)=R(\alpha,t)L(t), \quad \| R\|_{L^1}=1,
\end{equation} 
where $\| \cdot \|_{L^1}$ is the standard $L^1$ norm for absolutely integrable functions. Thus, the reparametrization problem here is reduced to finding the non-dimensional factor $R$, and this fact allows us to regard a monitor function $\varphi$ as complexity of a planar curve of the length $2\pi$ by the change of variables
\begin{equation}
\label{eq:def_change_s}
s=({L}/{2\pi})s', \quad s' \in [0,2\pi].
\end{equation}
Moreover, noting that multiplying $\varphi$ by a constant equally affects both sides of the definition (\ref{eq:def_equidist}), the $L^1$-type normalization
\begin{equation}
\label{eq:def_normalize}
\varphi^* = 2\pi({\varphi}/\| \varphi \|_{L^1}),
\end{equation}
is applied after the rescaling (\ref{eq:def_change_s}). As we see below, this modification, which is apparently meaningless in theory, periodizes derivatives of the arclength parametrization of the planar curve with the curvature $\kappa \equiv \varphi^*$ and plays a crucial role in spatial discritization in Section \ref{method}.\par
Now consider the interpolated curvature
  \begin{equation}
  \label{eq:def_dummy}
    \kappa(s,t) =(1-t)\kappa_0(s)+ t\kappa_1(s), \quad s\in[0,2\pi],\quad t\in[0,1],
  \end{equation}
where the two end-points are $\kappa_0 \equiv 1$ and $\kappa_1 \equiv \varphi^*$, respectively. By assuming that $L(t) = 2\pi$ for all $t\in[0,1]$, the linear interpolation  (\ref{eq:def_dummy}) defines a one-parameter family $\{\mathbf{X}(t)\}$ of planar curves up to a rigid transformation, or equivalently ``time evolution" of a single curve where we construct the desired function $R$. In general, motion of the generated curve is nonlinear, but controlling its shape at the curvature level avoids forming a cusp that we may encounter if coordinates $(x,y)$ are directly interpolated. Besides, in general each $\mathbf{X}(t)$ is open for $t>0$, while the initial state $\mathbf{X}(0)$ is the unit circle. However, note that the end-point curvatures satisfy
  \begin{equation}
  \label{eq:end_property}
    \kappa_0,\,\kappa_1>0,\quad \int_0^{2\pi} \kappa_0 ds = \int_0^{2\pi} \kappa_1 ds=2\pi,
  \end{equation}
and therefore the same properties hold for all $t\in[0,1]$. Since the angle $\theta$ is obtained by integrating the relation $\theta_s=\kappa$, we have
    \begin{equation}
      \label{eq:def_angle_int}
    \theta (t) =\partial^{-1}_s[ \kappa(t)],\quad \theta_t = \partial^{-1}_s [(\kappa_1-\kappa_0)(t)],
    \end{equation}
and the components of the tangent vector $\mathbf{X}_s$ are 
    \begin{equation}
     \label{eq:def_tvecor_dummy}
     x_s(t) =\cos\theta(t),\quad y_s(t) =\sin\theta(t).
  \end{equation}
Thus, the normalization (\ref{eq:def_normalize}) guarantees the periodicity $\mathbf{X}_s(0,t)=\mathbf{X}_s(2\pi,t)$ for all $t$ and typically reduces complexity of the final state $\mathbf{X}(1)$. Setting the starting point $\mathbf{X}(0,t)=(0,0)$ for all $t$, each coordinate is obtained by the anti-derivatives
\begin{equation}
     x(t) =\partial^{-1}_s [\cos\theta(t)],\quad y(t) =\partial^{-1}_s [\sin\theta(t)],
\end{equation}
and hence the components of the velocity $\mathbf{X}_t$ are
\begin{equation}
     x_t(t) =-\partial^{-1}_s [\theta_t \sin\theta(t)],\quad y_t(t) =\partial^{-1}_s [\theta_t \cos\theta(t)].
\end{equation}
To solve the evolution equations (\ref{eq:evol_theta_sa}), we need information on the normal velocity $U$ and the tangential velocity $V$ satisfying (\ref{eq:diff_Vphi}).  To this end, taking the dot product $\mathbf{X}_t \cdot \mathbf{n}$, we obtain
          \begin{align}
\label{eq:dummy_U}    U(t) &=-y_s \cdot x_t + x_s \cdot y_t, 
  \end{align}
  and differentiating both sides of (\ref{eq:dummy_U}) with respect to $s$ leads to
            \begin{align}
\label{eq:dummy_Us}     U_s(t) =-\kappa(x_s x_t +y_s y_t  ) + \theta_t.
  \end{align}
Moreover, the second derivative $U_{ss}$ is necessary for representing the term $V_\alpha$:
            \begin{align}
\label{eq:dummy_Uss}    U_{ss}(t) 
=-\kappa^2 U(t) +\frac{\kappa_s}{\kappa}\{U_s(t) - \theta_t\}+\kappa_t.
  \end{align}
Having these velocity data, the tangential component $V$ and its derivative $V_\alpha$ are computed for the equidistribution rule (\ref{eq:def_equidist}) with the special case $\varphi(\kappa) = \kappa$. First, note that Eq.\,(\ref{eq:def_diff_sa_phi}) simplifies as 
 \begin{equation}
\label{eq:def_diff_sa_phi_simple}
(s_\alpha \varphi)_t = (\kappa V+U_s)_\alpha.
  \end{equation}
Here, since the planar curve $\mathbf{X}(t)$ has the moving end-point $\mathbf{X}(2\pi,t)$ while $\mathbf{X}(0,t)$ is fixed at $(0,0)$, the boundary values of the tangential velocity $V$ are $V(0,t)=0$ and
\begin{align}
\nonumber V(2\pi,t)&=x_s x_t +y_s y_t  \\
\label{eq:def_Vend}             &=-\frac{1}{\kappa}(U_s(2\pi,t)-\theta_t(2\pi,t))\\
\nonumber             &=-\frac{1}{\kappa}U_s(2\pi,t),
\end{align}
where we used (\ref{eq:dummy_Us}) and then (\ref{eq:end_property}), and the integral of $(s_\alpha \varphi)_t$ over $[0,2\pi]$ is therefore
 \begin{equation}
 \int_0^{2\pi}(s_\alpha \varphi)_t d\alpha = \int_0^{2\pi} (\kappa V+U_s)_\alpha d\alpha =0.
  \end{equation}
Thus, from (\ref{eq:def_equidist}) and (\ref{eq:def_diff_sa_phi_simple}), the equation for the tangential velocity $V$ is reduced to
 \begin{equation}
(\kappa V+U_s)_\alpha=0,
  \end{equation}
and we obtain its explicit representation in terms of $U_s$ and $\kappa$
 \begin{equation}
V = -\frac{1}{\kappa}\partial^{-1}_\alpha [(U_s)_\alpha]  = -\frac{1}{\kappa}U_s(\alpha,t),
  \end{equation}
which is consistent with the boundary conditions as $\mathbf{X}_t(0,t)=(0,0)$ and $\theta_t(0,t)=0$ in (\ref{eq:dummy_Us}), and its first derivative with respect to $\alpha$
\begin{equation}
V_\alpha =-\frac{\kappa_\alpha}{\kappa}V -\frac{1}{\kappa}(U_s )_\alpha =-s_\alpha \frac{\kappa_s}{\kappa}V -\frac{s_\alpha}{\kappa}U_{ss}.
  \end{equation}
As opposed to (\ref{eq:diff_Vphi}) and (\ref{eq:int_V}), note here that the formulae for $V$ and $V_\alpha$ with $\varphi \equiv \kappa$ are purely algebraic provided the other data are given. This property is crucial in practical implementations because it is hard to efficiently approximate integrals involving the non-periodic factor $U$ and its derivatives.


\section{Numerical methods}\label{method}
This section is devoted to describing numerical methods that solve the reparametrization problem in Section \ref{formulate}. We first introduce building blocks of spectrally accurate calculus and its acceleration, and then develop the main algorithm that harnesses advantages of the interpolation-based approach.
\subsection{Spatial discretization}
In our computations, functions on the interval $[0,2\pi]$ are sampled at equispaced $N$ points
\begin{equation}
\label{eq:def_grid}
\alpha_j = jh,\quad h=\frac{2\pi}{N}, \quad j=0,1,\ldots, N-1.
\end{equation}
Then, the Fourier coefficients $\hat{f}$ of a smooth $2\pi$-periodic function $f$ are approximated by the trapezoidal rule with spectral accuracy:
\begin{equation}
\label{eq:def_ft_alpha}
\hat{f}(k)=\frac{1}{2\pi}\int_0^{2\pi}f(\alpha)e^{- i k{\alpha}}d\alpha \approx \frac{1}{N} \sum_{j=0}^{N-1} f(\alpha_j) e^{- i k{\alpha_j}}, 
\end{equation}
and the standard fast Fourier transform (FFT) is effective for summing up (\ref{eq:def_ft_alpha}) up to the Nyquist frequency. With these data, the derivative of $f$ is obtained by differentiating its truncated Fourier series term by term:
\begin{equation}
\label{eq:def_spec_diff}
f_\alpha(\alpha)\approx  \sum_{k=-N/2}^{N/2} (ik) \hat{f}(k) e^{ i k{\alpha}}.
\end{equation}
Similarly,  the anti-derivative $\partial^{-1}_\alpha f$ is approximated by
\begin{equation}
\label{eq:def_spec_int}
\int_0^\alpha f' d\alpha'\approx \alpha \hat{f}(0) + \sum_{\substack{k\neq 0 \\ k=-N/2}}^{N/2-1} \frac{\hat{f}(k)}{(ik)} (e^{ i k{\alpha}}-e^{0}),
\end{equation}
where the exponential sum is computed by the inverse FFT with $\hat{f}(0)=0$. In the following, the formulae (\ref{eq:def_spec_diff}) and (\ref{eq:def_spec_int}) are called spectral differentiation and integration, respectively. Here, it should be noted that performing calculus in the Fourier space requires just $\mathcal{O}(N)$ operations and the main computational efforts are devoted to the summations in the forward and inverse transformations.
\par
On the other hand, we also consider the Fourier coefficients of periodic functions in the arclength parametrization
\begin{equation}
\label{eq:def_ft_arc}
\mathcal{F}[f](k)=\frac{1}{L}\int_0^{L}f(s)e^{-2\pi i k\frac{s}{L}}ds,
\end{equation}
where $L$ denotes the total arclength of a periodic curve $\mathbf{X}$. These quantities can be accessed via the natural change of variables $ds = s_\alpha d\alpha$:
\begin{equation}
\label{eq:changevar}
\frac{1}{L}\int_0^{L}f(s)e^{-2\pi i k\frac{s}{L}}ds = \frac{1}{L}\int_0^{2\pi}(f(s(\alpha))s_\alpha)e^{-2\pi ik \frac{s(\alpha)}{L}}d\alpha,
\end{equation} 
and the trapezoidal rule with the data on (\ref{eq:def_grid}) leads to the formula
\begin{equation}
\label{eq:def_dft_arc}
\mathcal{F}[f](k)\approx \frac{h}{L} \sum_{j=0}^{N-1} (f(s(\alpha_j))s_\alpha (\alpha_j)) e^{-2\pi i \frac{s(\alpha_j)}{L}}.
\end{equation}
Here, it should be noted that nonuniform sample points $\{s(\alpha_j)\}$ may resolve functions of the arclength successfully, whereas complex exponentials, which oscillate uniformly in $[0,L]$, can be underresolved for sufficiently large $k$. To resolve both factors, functions of the variable $\alpha$ are interpolated by the truncated Fourier series to an upsampled grid with $N_\text{up}$ points in $[0,2\pi]$ prior to computing the coefficients $\mathcal{F}[f]$. We remark here that this interpolation only improves resolution in relation to the complex exponentials, because temporal data between the grid points (\ref{eq:def_grid}) are given not from the analytical data but from the initial $N$-point representation.\par
Unlike the formula (\ref{eq:def_ft_alpha}), the discrete sum (\ref{eq:def_dft_arc}) is not compatible with the standard FFT, because in general the nodes $\{s(\alpha_j)\}$ are distributed unevenly in $[0,L]$. Instead, an efficient algorithm called the {\it nonuniform fast Fourier transform} (NUFFT) is capable of approximating the following types of sums to a prescribed relative accuracy $\epsilon_\text{rel}$: 
\begin{alignat}{2}
\label{eq:def_nufft_type1}&\mbox{(Type-1)}\quad &\hat{f}_k=\sum_{j=0}^{N-1}f_je^{-ik x_j}, \quad &k=-\frac{M}{2},\ldots, \frac{M}{2}-1, \\
\label{eq:def_nufft_type2}&\mbox{(Type-2)}\quad &f_j=\sum_{k=-\frac{M}{2}}^{\frac{M}{2}-1}{\hat{f}_k}e^{ik x_j}, \quad &j=0,1,\ldots, N-1.
\end{alignat}
For instance, the Type-1 NUFFT regards the exponential sum (\ref{eq:def_nufft_type1}) as the exact Fourier transform of a measure $\mu$ of weighted Dirac masses, and replaces it with
\begin{equation*}
 \mu * \psi = \sum_{j=0}^{N-1}f_j \psi(\cdot-x_j), \quad \psi: \mbox{Gaussian, Kaiser-Bessel, etc.},
\end{equation*}
where $*$ is the convolution on the torus $\mathbb{R}/2\pi\mathbb{Z}$ and $\psi$ is a smooth periodic function with a compact numerical support (e.g., Gaussian and Kaiser-Bessel kernels). Then, the Fourier transform of the regularized measure is computed by the trapezoidal rule with the FFT acceleration, and dividing the results by the symbols $\hat{\psi}(k)$ yields approximations of the desired sums. For the Type-2 NUFFT, one may reverse the process above to approximate the sums (\ref{eq:def_nufft_type2}). Here, it is easy to see that the formula (\ref{eq:def_dft_arc}) naturally fits in the framework of the Type-1 NUFFT. On the other hand, once the Fourier coefficients of a function of the arclength are obtained, it can be interpolated back to the nonuniform grid $\{s(\alpha_j)\}$ using the Type-2 NUFFT. The convergence of such algorithms is first analyzed by Dutt and Rokhlin \cite{DuRo1993} for the Gaussian kernel, and now there are several general-purpose implementations including CMCL NUFFT \cite{GrLe2004}, NFFT3 \cite{KeKuPo2009}, and FINUFFT \cite{BaMaKl2019}. Among those, FINUFFT, which chooses $\psi$ to be the ``exponential of semicircle" kernel with nearly optimal aliasing errors \cite{Barnett2021}, is currently the fastest library without precomputation, and it is employed in our code that repeatedly invokes the Type-2 NUFFT with moving nodes $\{s(\alpha_j)\}$.

\subsection{Main algorithm} Now we assemble the above components into a single algorithm that perform reperametrization with the interpolation (\ref{eq:def_dummy}). The input data here are a periodic planar curve $\mathbf{\Gamma}$ with $s_\alpha>0$ and a monitor function $\varphi$, which are both assumed to be smooth periodic functions of the variables $\alpha$ and $s$.\par 
\textbf{Step 1}: The algorithm extracts a canonical parametrization from the curve $\mathbf{\Gamma}$ and stores it as intermediate data before switching to a refined representation. Although there can be arbitrariness for this step, the arclength parametrization (\ref{eq:X_arclength}) is considered the most natural one, because it is linked to the metric of a planar curve and therefore cannot be lost by a change of parametrizations \cite{Koga2021}. An immediate consequence of this invariance is the change of variables (\ref{eq:changevar}), which allows direct access to  functions of the arclength variable from arbitrary parametrization. Firstly, the derivatives of the coordinates $(x,y)$ are approximated by spectral differentiation, and the local spacing $s_\alpha$ and the curvature $\kappa$ are computed from those data. Then, after obtaining the arclength $s$ by spectral integration, the Fourier coefficients $\mathcal{F}[\kappa]$ are discretized by the trapezoidal rule (\ref{eq:def_dft_arc}), and the resulting sums are approximated by the Type-1 NUFFT up to the wavenumber $k=N_1/2$, where $N_1$ is the sample size for the intermediate data. Finally, the arclength parametrization of $(x,y)$ is obtained by applying spectral integration to the relation $\theta_s=\kappa$ and (\ref{eq:def_theta}), and their Fourier coefficients are the desired intermediate data. This uniformalization at the curvature level is due to the fact that many sophisticated simulations of moving boundary problems are formulated with the variables $(\theta, s_\alpha$), which it is somewhat tricky to compute according to the direct formula $\theta = \arctan (y_\alpha/x_\alpha)$, as pointed out in \cite{HoLoSh1994}. \par
\textbf{Step 2}: A new local spacing $s_\alpha$ is constructed by solving the initial-value problem that arises from the curvature interpolation (\ref{eq:def_dummy}). Due to the normalization (\ref{eq:def_normalize}), the normal velocity $U$ and its derivatives are written in terms of periodic functions and their anti-derivatives, and hence spatial discretization in this simulation is exponentially convergent for spectral differentiation and integration. Here, we use $N_2$ points for discretizing the one-parameter family of curves $\{\mathbf{X}(t)\}$ and the solutions to  the system (\ref{eq:evol_theta_sa}), where the value of $N_2$ is typically smaller than $N_1$. On the other hand, since the angle $\theta(t)$ is given at every $t\in [0,1]$, the algorithm solves only the second equation in (\ref{eq:evol_theta_sa}), and the local spacing $s_\alpha$ is evolved in time using the classical fourth-order Runge-Kutta method. To obtain the known data in the time evolution, their truncated Fourier series are evaluated at nonuniform nodes $\{s(\alpha_j)\}$, and those operations are accelerated by the Type-2 NUFFT. After completing the simulation,  the non-dimensional factor $R$ in (\ref{eq:def_R}) is found and can be directly used for $\mathbf{\Gamma}$ because it is shared between the monitor function $\varphi$ and its rescaled counterpart $\varphi^*$.  \par
\textbf{Step 3}: We obtain a new parametrization of the curve $\mathbf{\Gamma}$ from the results of the previous two steps. The local spacing $s_\alpha$ from Step 2 gives rise to a new mapping $s:[0,2\pi]\rightarrow [0,L]$, and its numerical approximation is done by spectral integration. Then, the invariant Fourier coefficients of the coordinates $(x,y)$ from Step 1 are inverted by evaluating the corresponding truncated Fourier series with the target points $\{s(\alpha_j)\}$ of size $N_3$. At this step, the output size $N_3$ must be equal to or larger than $N_2$ in order that the new discrete representation resolves both $s_\alpha$ and the curve itself, while it may be  smaller than $N_1$ due to effects of mesh refinement. As above, the evaluations of the truncated Fourier series can be accelerated by the Type-II NUFFT.    \par 
As a validation of the algorithm, improvements to spatial resolution are evaluated by recomputing the Fourier coefficients of the arclength parametrization from downsampled data and comparing those against the intermediate data at Step 1. That is, if spatial resolution of the new parametrization is superior to the input, it should be possible to compute the invariant Fourier coefficients with a reduced number of discrete points, which implies faster global convergence to the original shape. Such a quantitative comparison is essential for claiming that the new discrete representation is in fact a refined approximation of the input curve, because the interpolated values from Step 3 are trivially correct pointwise but not necessarily as a whole.\par 
Here, it should be mentioned that many alternatives have been proposed for computing the arclength parametrization from a given parametric curve. Examples include iterative schemes based on the Newton's method \cite{BeRo2014, BaSh1990, HoLoSh1994} or variants of de Boor's algorithm \cite{MaNoRoIn2019, SeLa2020}, which search for a set of values $\{\alpha_j\}$ corresponding to $\{s(\alpha_j)\}$ equispaced in $[0,L]$, and asymptotic methods \cite{SeYa2011, MiSe2004}, which slides discrete points along a planar curve by $\mathbf{X}_t=V\mathbf{t}$ until the local spacing $s_\alpha$ converges up to a given tolerance. Nevertheless, we prefer the non-iterative approach \cite{Koga2021} based on the coefficients (\ref{eq:def_dft_arc}) because its performance is determined solely by resolution properties of the trapezoidal rule, whereas additional errors from numerical optimization or temporal discretization are inevitable for the others.


\section{Numerical results}\label{results}
Now we show numerical results to validate the algorithm described in Section \ref{method}. In the following, the implemented codes are run on MacBook Pro 2020 with Intel Core i7-1068NG7 2.3Ghz and 16GB memory, and the relative accuracy for the NUFFT algorithms is fixed to $\epsilon_\text{rel}=10^{-15}$. 
%
%
%
%
%
%
%
\subsection{Convergence study}
\begin{figure}[t]
\begin{center}
 \includegraphics[width=0.8\linewidth]{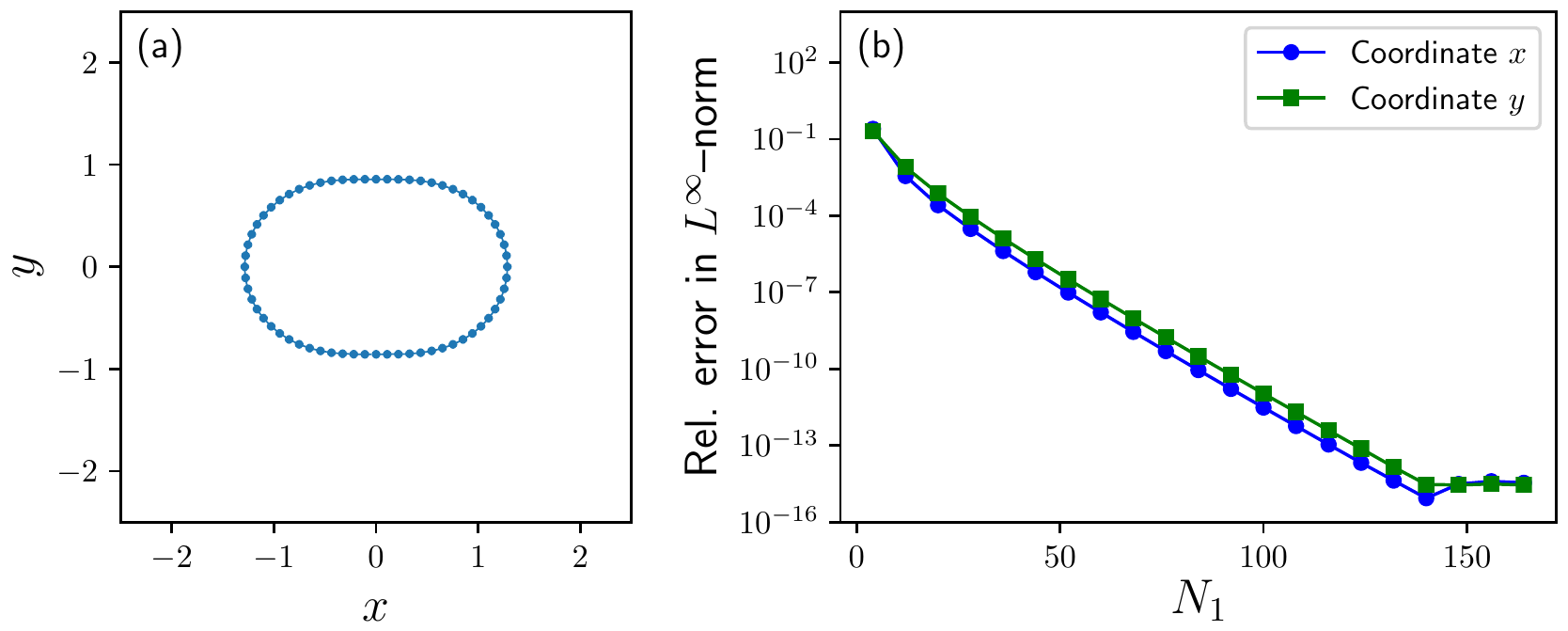}
\end{center}
\caption{\label{fig:fig_conv_canon} Example (\ref{eq:linear_period}) and its spatial discretization: (a) representation in $(x,y)$-plane, and (b) spectral accuracy in computing arclength parametrization at Step 1.} 
\end{figure}

We start numerical experiments by checking spectral accuracy of Step 1. This is essentially a reproduction of the convergence study in \cite{Koga2021}, where the technique is tested on the following example in polar coordinates $(r,\eta)$:
%
\begin{equation}
\label{eq:linear_period}
r(\eta)=1+\epsilon_P P_2\biggl(\cos\biggl(\eta-\frac{\pi}{4}\biggr)\biggr), \;\; \eta \in [0,2\pi].
\end{equation} 
Here, $\eta$ is the polar angle, $r$ is the radius as a function of $\eta$, and $P_2$ is the Legendre polynomial of order two. This specific choice corresponds to the zero-velocity state of  periodic motion of the linearized droplet dynamics \cite{LuMa1988}. The representation in the polar coordinates allows to compute the inverse mapping $(x,y)\mapsto \eta$ and evaluate accuracy of numerical results against the analytical data (\ref{eq:linear_period}). In Fig.\,\ref{fig:fig_conv_canon}, the planar curve (\ref{eq:linear_period}) and the convergence of the scheme are shown for $\epsilon_P=2/7$ and $N_\text{up}=N_1$. 
%
As shown in Fig.\,\ref{fig:fig_conv_canon}(a), it has the variable radius $r$ and therefore the local spacing  $s_\eta=(r^2+r_\eta^2)^\frac{1}{2}$ is not constant. Next, in Fig.\,\ref{fig:fig_conv_canon}(b), numerical errors in $x$ and $y$ are plotted versus $N_1$ in terms of the $L^\infty$ norm divided by the maximal value of each function. The exponential convergence of the trapezoidal rule is evident, and the errors reache a level slightly above $10^{-15}$ for $N_1\geq 140$.\par
\begin{figure}[t]
\begin{center}
 \includegraphics[width=0.8\linewidth]{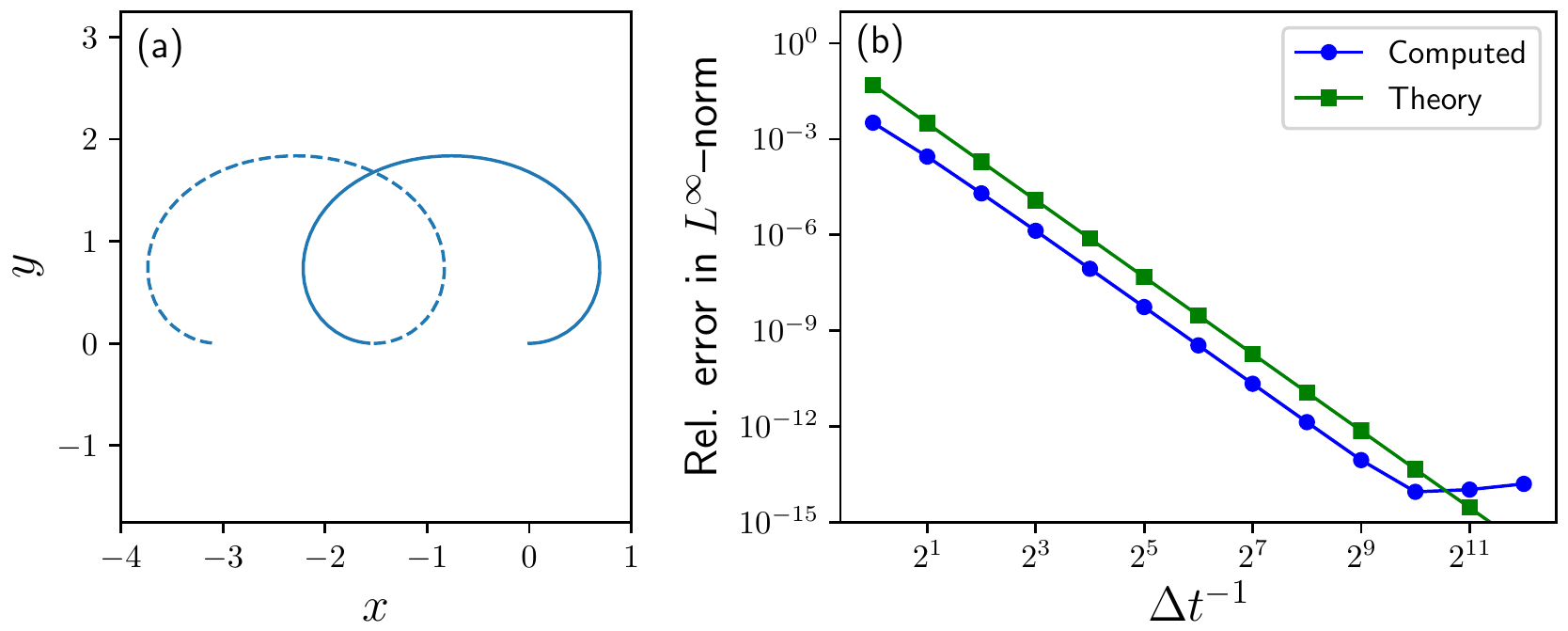}
\end{center}
\caption{\label{fig:fig_dummy_conv} Curvature interpolation and its time discretization: (a) end-point curve corresponding to (\ref{eq:monitor_simple}) after normalization (\ref{eq:def_normalize}), and (b) 4th-order accuracy of Runge-Kutta method at Step 2.} 
\end{figure}

Next, we examine accuracy of spatial and temporal discretization in numerical simulations of Step 2 and claim effectiveness of the normalization (\ref{eq:def_normalize}). For this purpose,  the simple monitor function
\begin{equation}
\label{eq:monitor_simple}
\varphi_0(s)=0.5 + 0.25\cos(s),\quad s\in [0,2\pi],
\end{equation}
is considered here as an intuitive example that highlights the core ideas behind the algorithm. Clearly, the function $\varphi$ is strictly positive, smooth, and periodic in the variable $s$, while its integral over $[0,2\pi]$ is just $\pi$ rather than $2\pi$. This means that the identification $\kappa_1\equiv \varphi_0$ yields the arclength parametrization of an open planar curve whose tangent vector $\mathbf{X}_s (1)$ is non-periodic in $s$. Then, the $L^1$-type nomralization (\ref{eq:def_normalize}) multiplies the function $\varphi$ by the factor $2\pi/\| \varphi \|_{L^1}$ so that $\mathbf{X}_s (1)$ rotates exactly once as $s$ moves from $0$ to $2\pi$. Figure\,\ref{fig:fig_dummy_conv}(a) shows the end-point curve $\mathbf{X}(1)$ generated from $\kappa_1 \equiv \varphi_0^{*}$ for the first two periods (solid line for the first and dashed line for the second). As seen, this object is not even closed and the point $\mathbf{X}(s,1)$ goes to infinity as $|s|\rightarrow \infty$. However, adjusting the monitor function results in periodic derivatives in the arclength parametrization and therefore is expected to give rise to spectral accuracy in computing velocity data of the one-parameter family $\{\mathbf{X}(t) \}$. Figure\,\ref{fig:fig_dummy_conv}(b) plots relative $L^{\infty}$ errors from simulations of the initial-value problem (denoted by ``Computed")  with the classical fourth-order Runge-Kutta method and $N_2=128$, and also shows a theoretical $\mathcal{O}(\Delta t^4)$ line as a reference (denoted by ``Theory"). In this convergence study, we use the right-hand of the definition (\ref{eq:def_equidist}) as the true solution, which is easily calculated by hand for the case (\ref{eq:monitor_simple}), and compare numerical approximations to the exact value. One can see that the error curve from the experiment is almost parallel to the $\mathcal{O}(\Delta t^4)$ line and ``saturated" at around the level of $10^{-14}$ with the rather small $N_2$. Thus, we conclude from these data that the normalization (\ref{eq:def_normalize}) enables spatial discretization of spectral accuracy in evolving the interpolated curve $\{\mathbf{X}(t) \}$, even though each $\mathbf{X}(t)$ is not periodic or closed for $t>0$.
%
%
%
%
%
%
\begin{figure}[t]
\begin{center}
 \includegraphics[width=\linewidth]{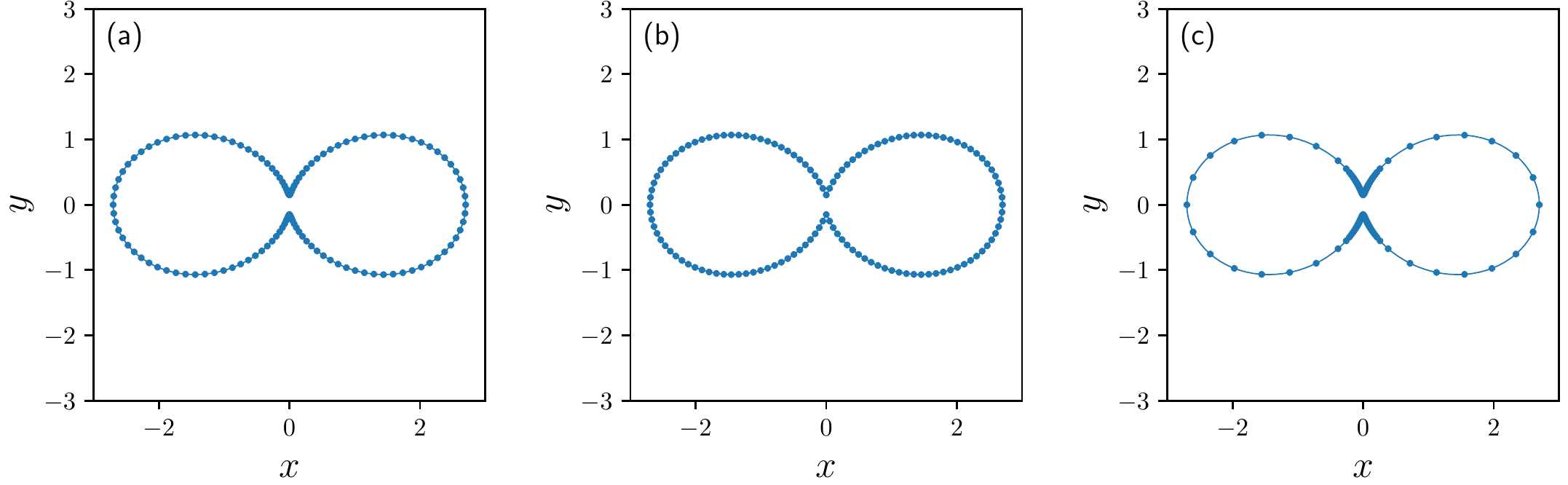}
\end{center}
\caption{\label{fig:fig_drop} Pinched droplet with different meshes: (a) input, (b) arclength, and (c) refined. } 
\end{figure}
\subsection{Pinched droplet}
Our first practical example arises from changing the value of $\epsilon_P$ in the analytical data (\ref{eq:linear_period}). That is, increasing the size of the deformation $\epsilon_P$ ``pinches the droplet" and seems to form two cusps as $\epsilon_P \rightarrow 2$. However, since the algorithm requires the input data to be smooth, we consider a value $\epsilon_P<2$ that corresponds to a nearly singular curve with locally high curvature. Figure \ref{fig:fig_drop}(a) shows the input curve with $\epsilon_P=1.7$ for which the maximal curvature $\kappa_\text{max}$ is approximately $220$. As seen in the plots, the original mesh has some concentration in the neighborhoods of the high-curvature regions and is relatively uniform otherwise. Such a distribution is due to the use of polar coordinates, where the local spacing $s_\eta$ is given by $s_\eta=(r^2+r_\eta^2)^\frac{1}{2}$ and therefore the points of the minimum radius often attain the minimum spacing. Following the recipe in Section \ref{method}, we first compute the Fourier coefficients (\ref{eq:def_ft_arc}) for the coordinates $(x,y)$ and store them as intermediate data, and its discretization in the the $x$-$y$ plane is shown in Fig.\,\ref{fig:fig_drop}(b). Here, those invariants are computed up to $k=16384$, and we observe that the sample points $N_1=32768$ is enough for discretizing the case $\epsilon_P=1.7$ and getting the relative $L^\infty$ error of $10^{-11}$  at Step 1. Unlike the case $\epsilon_P = 2/7$, it is necessary to set $N_{\text{up}}=2N_1$ for $\epsilon_P = 1.7$, because $N_{\text{up}}=N_1$  leads to poor accuracy in computing the Fourier coefficients of the arclength parametrization whose integrands have uniformly oscillating factors.\par
\begin{figure}[t]
\begin{center}
 \includegraphics[width=0.8\linewidth]{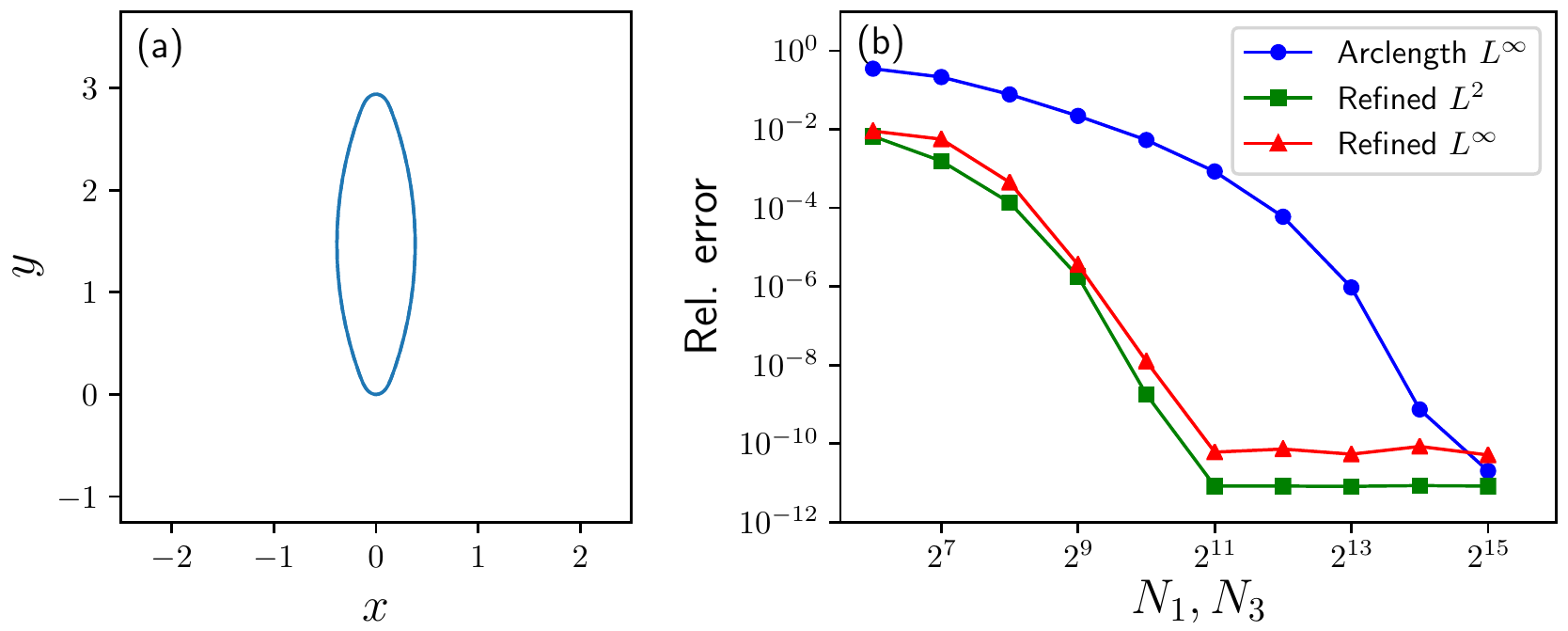}
\end{center}
\caption{\label{fig:fig_refine_conv_pinch} Improvements to spatial resolution for pinched droplet: (a) end-point curve corresponding to (\ref{eq:def_monitor_drop}), and (b) error plots showing exponential convergences with different slopes.} 
\end{figure}

Next, we design a suitable monitor function that resembles complexity of the input curve. Obviously, slowly oscillating functions such as  (\ref{eq:monitor_simple}) does not fit in the current situation, because the high-curvature regions of the curve are localized in the neighborhoods around the two isolated points. To this end, the present work writes practical monitor functions as superpositions of periodized Gaussian kernels whose locality is easily controlled by changing coefficients in the exponents. For the case $\epsilon_P=1.7$, the monitor function $\varphi_1$ is defined as
\begin{equation}
\label{eq:def_monitor_drop}
\varphi_1(s)= 1 +37\sum_{j\in \mathbb{Z}}e^{-\{7.5(s+2\pi j)\}^2}+37\sum_{j\in \mathbb{Z}}e^{-\{7.5(s-\pi+2\pi j)\}^2},
\end{equation}
where parameters are sought by trials and errors so that the later validation gives satisfactory results. In our implementations, the summation in (\ref{eq:def_monitor_drop}) is truncated when contributions from periodic images are below $10^{-15}$. The corresponding end-point curve $\mathbf{X}(1)$ is plotted in Fig.\,\ref{fig:fig_refine_conv_pinch}(a). Due to the symmetry of the data (\ref{eq:linear_period}), the resulted object turns out to be closed and enclose a convex domain, and the normalization (\ref{eq:def_normalize}) reduces its complexity for the case $\varphi_1$ with $\| \varphi_1 \|_{L^1} \approx 23$. For reference, a single-thread simulation at Step 2, where we choose $N_2=2048$ and $\Delta t= 10^{-4}$, takes approximately 2.4 minutes to achieve the relative $L^\infty$ error of $5\times 10^{-13}$ in the product $s_\alpha\varphi_1$. We do not observe speed gain from increasing the number of threads, which is because the order of $N_2$ is small in comparison to the main target of the NUFFT algorithms with OpenMP (see \cite{BaMaKl2019} for details). In fact, the computation time for the case $\kappa_0 \equiv 1$ is not significant in practical simulations of moving boundary problems with a time-dependent monitor function $\varphi$. Namely, once the equidistribution rule is well approximated for $\varphi(t)$, it is far less expensive to find the desired parametrization for $\varphi(t+\Delta t)$ by replacing the end-point curvatures $(\kappa_0, \kappa_1)$ in the interpolation (\ref{eq:def_dummy}) with the pair $(\varphi^*(t), \varphi^*(t+\Delta t))$. \par
Step 3 outputs a refined representation of the input data obtained by interpolating the arclength parametrization onto the nonunfirom grids $\{s(\alpha_j)\}$ that corresponds to the equidistribution rule with (\ref{eq:def_monitor_drop}). Figure \ref{fig:fig_drop}(c) shows discretization of the resulted parametrization with the same number of points as in Fig.\,\ref{fig:fig_drop}(a). It is evident that the high-curvature regions are more densely resolved and the mesh distribution is sparse away from the targets. However, it is not possible to conclude from this visualization that the refined parametrization can represent the planar curve with a smaller number of sample points than the original form (\ref{eq:linear_period}), and more quantitative comparisons between different meshes are necessary.\par
 As described in Section \ref{method}, we confirm such improvements by recomputing the Fourier coefficients (\ref{eq:def_ft_arc}) from downsampled data with varying $N_3$ and compare those against the intermediate data of full accuracy. Thus, the following quantitative analysis is essentially done against the arclength parametrization rather than the input data, which typically results in the worst among the three. In Fig.\,\ref{fig:fig_refine_conv_pinch}(b), three error curves are plotted on the logarithmic scale against the sample sizes $N_1$ and $N_3$. The blue curve (denoted by ``Arclength $L^\infty$") is for relative $L^\infty$ errors in the coordinates $(x,y)$ from comparing the arclength parametrization to the analytical data after inverting the results from Step 1, and it shows exponential convergence to the input data as well as a rough estimation on the best possible accuracy that can be obtained after reparametrization. Again, we set $N_{\text{up}}=2N_1$ to avoid inaccuracy in computing the Fourier coefficients with the uniformly oscillating factors. As opposed to the case $\epsilon={2}/{7}$, a few digits of accuracy are lost at the level of saturation, which is presumably due to errors from spectral differentiation with large $N_{\text{up}}$. On the other hand, the green curve (denoted by ``Refined $L^2$") is for relative $L^2$ errors that are computed via the Plancherel theorem, which amounts to computing the $l^2$ norm of the difference between the the intermediate data and the recomputed Fourier coefficients. As clear from the figure, the errors from the refined representation decay faster than those from the arclength parametrization, and we claim that the full accuracy with $N_1=32768$ is obtained by the refined representation with $N_3= 2048$ and $N_{\text{up}} = 65536$. Here, we remark that an integral norm tends to give smaller errors than pointwise evaluations. Nevertheless, Figure \ref{fig:fig_refine_conv_pinch}(b) also shows the red curve (denoted by ``Refined $L^\infty$") for relative $L^{\infty}$ errors computed after inverting the invariant Fourier coefficients from the new parametrization, which leads us to the same conclusion as above.   \par
 \begin{figure}[t]
\begin{center}
 \includegraphics[width=0.8\linewidth]{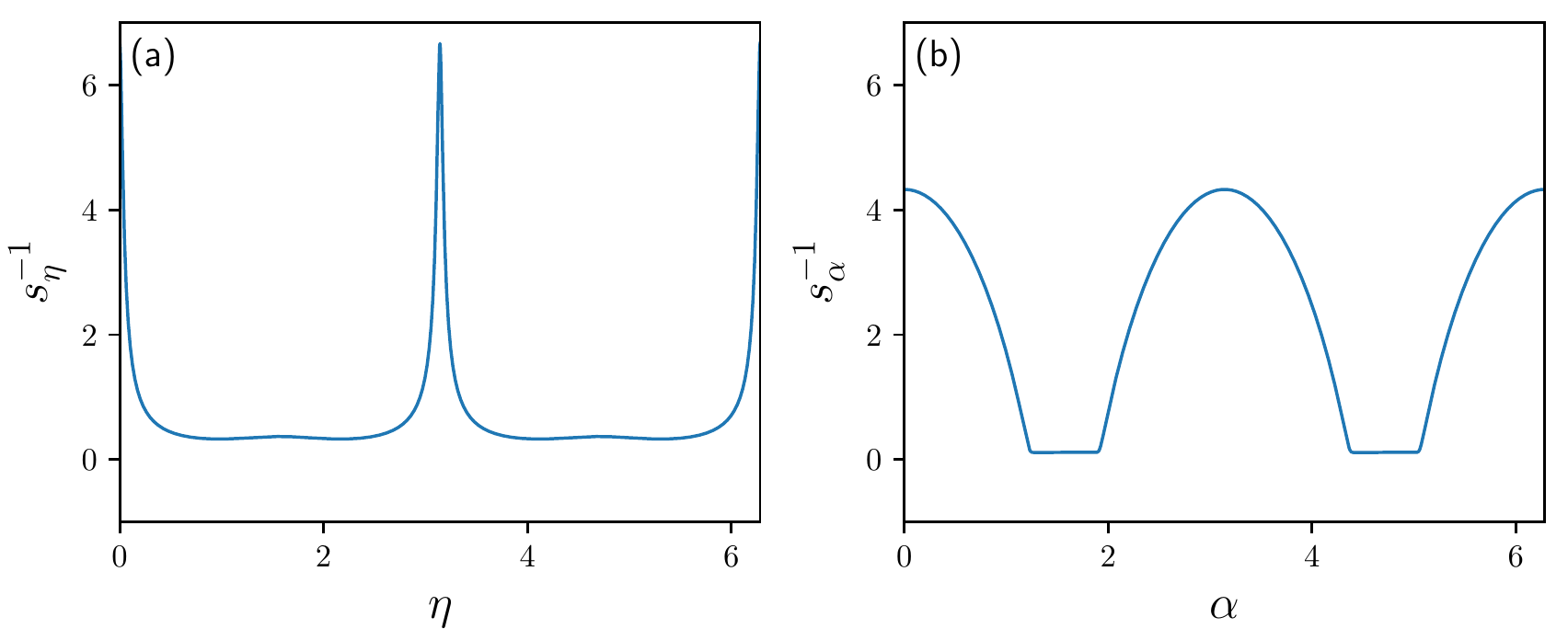}
\end{center}
\caption{\label{fig:fig_plot_space_pinch} Inverse of local spacing for pinched droplet as functions of different spatial variables: (a) polar angle $\eta$, and (b) equidistribution $\alpha$.} 
\end{figure}
Lastly, we confirm that the monitor function (\ref{eq:def_monitor_drop}) does not generate a pathological mesh distribution. For this purpose, Figure \ref{fig:fig_plot_space_pinch}(a) shows the reciprocal of the local spacing $s_\eta$ of the polar representation (\ref{eq:linear_period}) as a function of the angle $\eta$. As clearly seen, it has a sharp peak at $\eta=0$ and $\eta=\pi$, which means that the initial mesh barely resolves the small portion around the ``waist of the droplet" and results in the local spacing that contains high-frequency modes. As pointed out in Section \ref{method}, the major advantage of the approach based on the invariants (\ref{eq:def_dft_arc}) is that the success of the algorithm is solely determined by spatial resolution of the equispaced grids (\ref{eq:def_grid}) and the Type-1 NUFFT avoids a severe growth of computational costs from increasing $N_1$ and $N_\text{up}$, whereas optimization-based methods may additionally suffer from sharpness of the input data. On the other hand, Figure \ref{fig:fig_plot_space_pinch}(b) shows the same quantity $s_\alpha^{-1}$ of the new representation as a function of the variable $\alpha$.  Although it still has relatively high-frequency modes, which is due to the fast-decaying function used in the monitor function (\ref{eq:def_monitor_drop}), the sharpness for the polar representation is no longer found in the local spacing from Step 2. As confirmed above, such improvements to the parametrization appear as the faster global convergence in Fig.\,\ref{fig:fig_refine_conv_pinch}(b).
%
%
%
\subsection{Periodized peakons}
\begin{figure}[t]
\begin{center}
 \includegraphics[width=\linewidth]{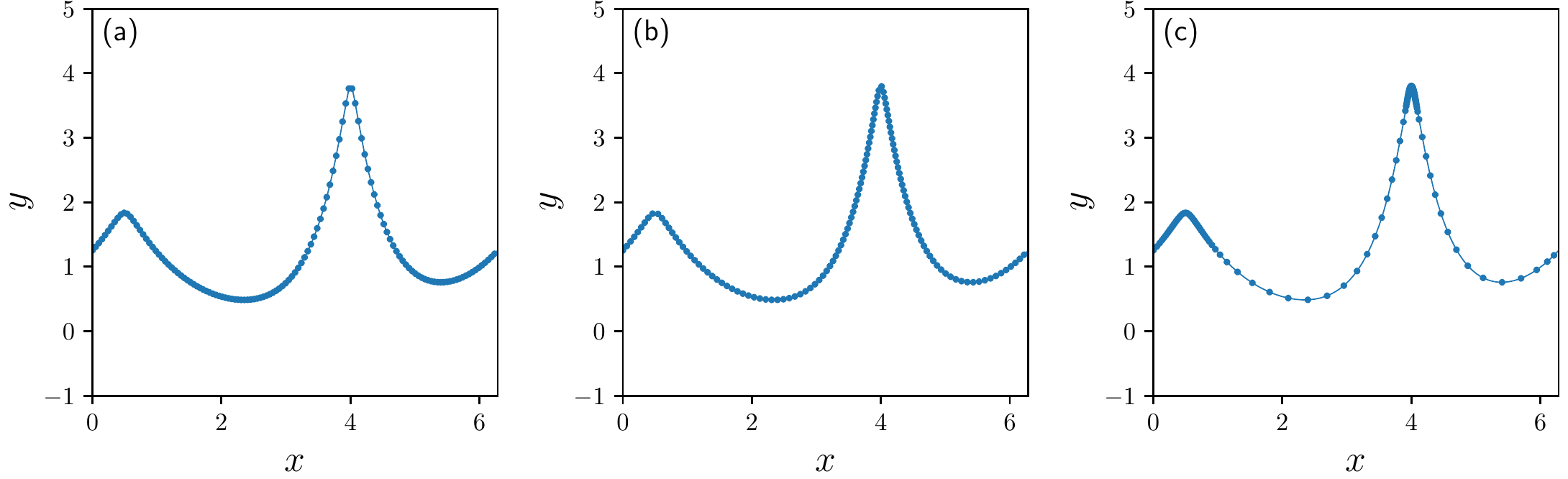}
\end{center}
\caption{\label{fig:fig_waves} Periodized peakons with different meshes: (a) input, (b) arclength, and (c) refined. } 
\end{figure}
The second example is the following superposition of periodic waves in the function representation $(x, y(x))$ on $[0,2\pi]$:
\begin{equation}
\label{eq:def_peakons}
y(x) =2\sum_{j\in \mathbb{Z}}e^{-\sqrt{(x-0.5+2\pi j)^2+\epsilon_\text{R}}} +4\sum_{j\in \mathbb{Z}} e^{-\sqrt{\{2(x-4+2\pi j)\}^2+\epsilon_\text{R}}}, \quad \epsilon_\text{R}>0.
\end{equation}
This specific wave form is found in the context of peakon solutions to the Camassa-Holm equation \cite{BeSaSz1999}, which are originally singular (i.e. $\epsilon_R=0$) and decay exponentially fast as $|x|\rightarrow \infty$. Hence, those functions are periodized here by the same technique as in (\ref{eq:def_monitor_drop}), while a larger number of periodic images are needed for the truncation with the threshold $10^{-15}$. The basic ideas of reparametrization applies to horizontally periodic curves without any modification, because the coordinate $x$ has the fixed boundary values as well as periodic derivatives and can be periodized by subtracting a linear function over $[0,2\pi]$. As opposed to the previous example, increasing the size of the rounding parameter $\epsilon_R$ smoothes the tip of each peakon and it forms a cusp when $\epsilon_R \rightarrow 0$, which requires us to consider strictly positive values of $\epsilon_R$ only. Figure \ref{fig:fig_waves}(a) shows the input curve (\ref{eq:def_peakons}) with $\epsilon_\text{R}=10^{-2}$ and the maximal curvature $\kappa_\text{max}\approx 144$. As seen in the plots, the initial mesh is rather sparse around the high-curvature regions. Such a distribution is due to the use of function representations, where the local spacing $s_x$ is given by $s_x= \sqrt{y_x^2+1}$ and therefore it is directly linked to the magnitudes of derivatives. In this ``ill-conditioned" case, it is found that Step 1 requires the upsampling $N_{\text{up}}=16N_1$ for discretizing the Fourier coefficients of the arclength parametrization with sufficient accuracy, whereas we choose $N_{\text{up}}=65536$ for $N_1 \geq 8192$ due to the limitation on the memory allocation. Again, those invariants are computed up to $k=16384$, and we observe that the sample points $N_1=32768$ is enough for the relative $L^\infty$ error slightly below $10^{-11}$ with $\epsilon_R=10^{-2}$. A discrete representation of the arclength parametrization is shown in Fig.\,\ref{fig:fig_waves}(b).\par
\begin{figure}[t]
\begin{center}
 \includegraphics[width=0.8\linewidth]{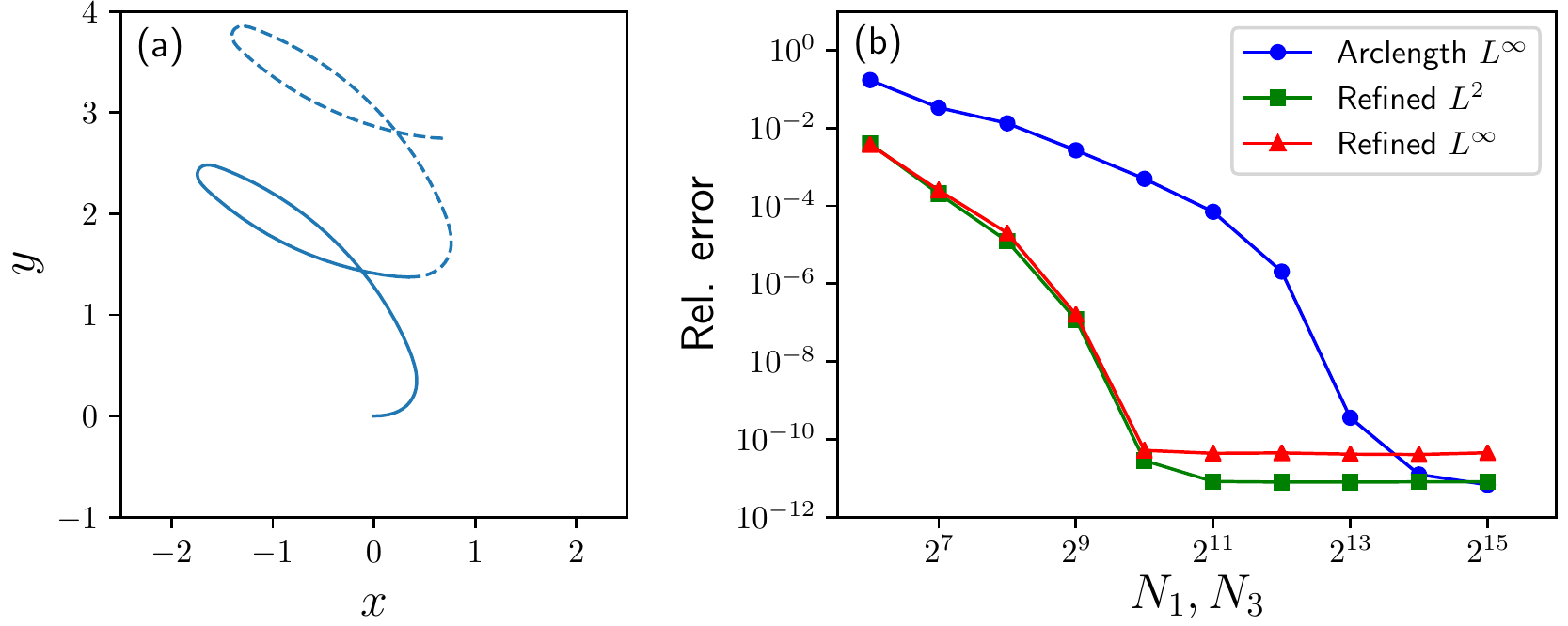}
\end{center}
\caption{\label{fig:fig_refine_conv_wave} Improvements to spatial resolution for periodized peakons: (a) end-point curve corresponding to (\ref{eq:def_monitor_waves}), and (b) error plots showing exponential convergences with different slopes. } 
\end{figure}

For the example (\ref{eq:def_peakons}), a monitor function $\varphi_2$ is also defined to be a superposition of periodized Gaussian kernels but with non-shared parameters:
\begin{equation}
\label{eq:def_monitor_waves}
\varphi_2(s)= 1 +10\sum_{j\in \mathbb{Z}}e^{-\{3(s-0.4+2\pi j)\}^2}+37\sum_{j\in \mathbb{Z}}e^{-[7.5\{s-(\pi+0.692)+2\pi j\}]^2}.
\end{equation}
Again, the constants in (\ref{eq:def_monitor_waves}) are sought by trials and errors in order to minimize errors at the validation step. In Fig.\,\ref{fig:fig_refine_conv_wave}(a), the resulting end-point curve $\mathbf{X}(1)$ is plotted for the first two periods. As opposed to (\ref{eq:def_monitor_drop}),  the absence of a symmetry leads to an open and self-intersecting curve, and the normalization (\ref{eq:def_normalize}) plays an essential role in allowing spatial discretization of spectral accuracy while reducing the complexity for the case $\varphi_2$ with $\| \varphi_2 \|_{L^1} \approx 21$. In this case, Step 2 uses the parameters $N_2=2048$ and $\Delta t= 5 \times 10^{-5}$ to achieve the relative $L^\infty$ error of $2\times 10^{-13}$ in the product $s_\alpha\varphi_2$.\par 
Now, we obtain a new representation of the input data by combining the outputs of the previous two steps. Figure \ref{fig:fig_waves}(c) shows discretization of the resulted parametrization with the same number of points as in Fig.\,\ref{fig:fig_waves}(a), and one can easily see that spatial resolutions around the two tips are significantly improved while the mesh distribution is sparse otherwise. Note here that the number of discrete points resolving each peakon can be adjusted by simply changing the parameters of the corresponding term in $\varphi_2$, which is the main advantage of choosing monitor functions to be linear combinations of smooth functions with compact numerical supports.  \par
 As before, improvements to spatial resolution are evaluated by reconstructing the arclength parametrization in the Fourier space from the refined reperesentation with variable $N_3$. In Fig.\,\ref{fig:fig_refine_conv_wave}(b), the blue curve (denoted by ``Arclength $L^\infty$") is for relative $L^\infty$ errors from comparing the arclength parametrization at Step 1 to the exact form (\ref{eq:def_peakons}), the green curve (denoted by ``Refined $L^2$") for relative $L^2$ errors that are computed via the Plancherel theorem, and the red curve (denoted by ``Refined $L^\infty$") for relative $L^{\infty}$ errors from the refined parametrization. Again, we set $N_{\text{up}}=16N_1$ for $N_1\leq 4096$ and $N_{\text{up}}=65536$ otherwise to avoid inaccuracy in computing the invariant Fourier coefficients. As one can see, the errors from the refined representation decay faster than those from the arclength parametrization, and we claim that the accuracy with $N_1=16384$ is obtained by the refined representation with $N_3= 2048$ and $N_{\text{up}}=32768$ in the $L^2$ sense.\par
\begin{figure}[t]
\begin{center}
 \includegraphics[width=0.8\linewidth]{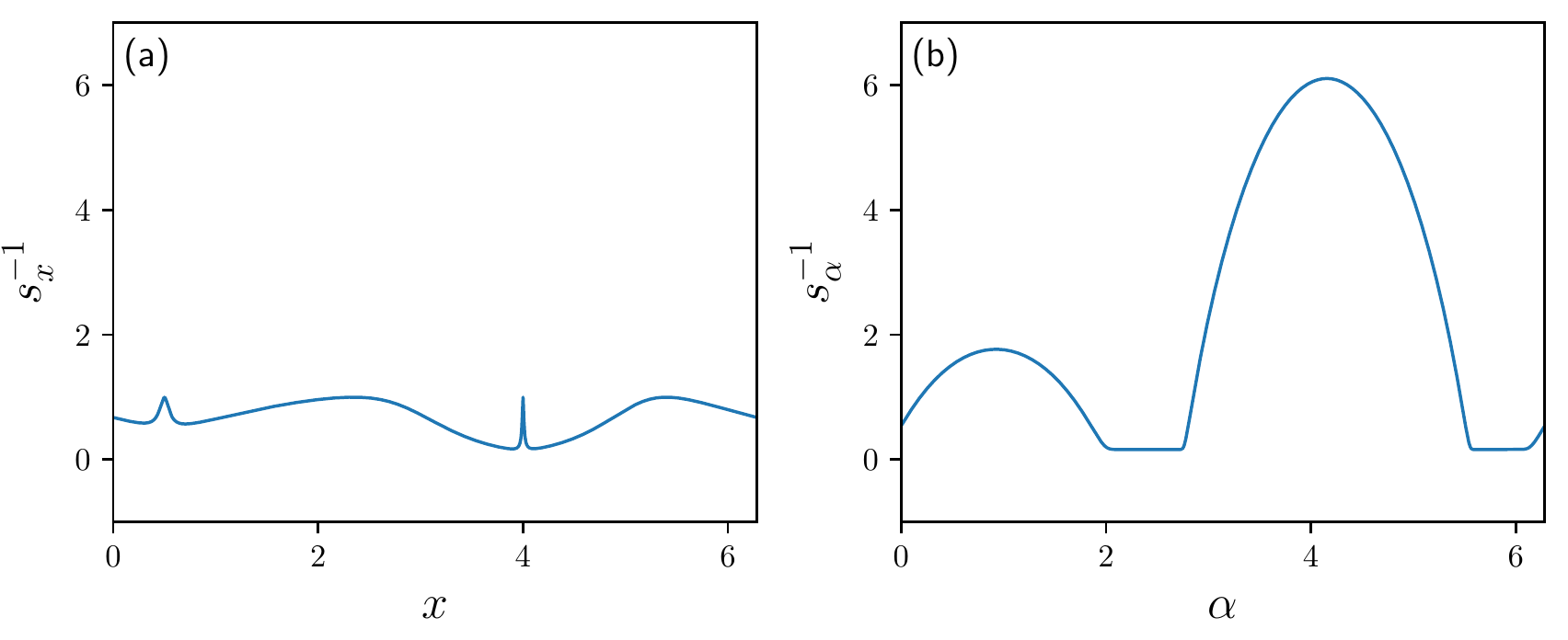}
\end{center}
\caption{\label{fig:fig_plot_space_wave} Inverse of local spacing for periodized peakons as functions of different spatial variables: (a) coordinate $x$, and (b) equidistribution $\alpha$.} 
\end{figure}

We finish the second example by comparing the local spacing $s_x$ of the input data to that of its refinement. Figure \ref{fig:fig_plot_space_wave}(a) shows the reciprocal of the local spacing $s_x$ of the representation (\ref{eq:def_peakons}) as a function of the coordinate $x$. Since the tip of each peakon is smoothed by the rounding parameter $\epsilon_R>0$, the derivative $y_x$ vanishes at $x=0.5$ and $x=4$ and hence the local spacing attains its local minima there. However, these minima are located within pathologically narrow intervals and small values around the points are barely sampled by the equispaced grids in $[0,2\pi]$, which results in the sparsity in the neighborhoods of the most complex geometry. As mentioned earlier, the invariant-based approach at Step 1 succeeds in converting the input curve (\ref{eq:def_peakons}) to its arclength parametrization. On the other hand, Figure \ref{fig:fig_plot_space_wave}(b) shows the same quantity $s_\alpha^{-1}$ of the new representation as a function of the variable $\alpha$.  As opposed to Fig.\,\ref{fig:fig_plot_space_wave}(a), the intervals with small values of the function $s_\alpha$ are wide enough and thus the uniform sampling in $\alpha$ resolves the high-curvature regions efficiently.

\section{Conclusion}\label{conclusion}
In this paper, we have developed a static algorithm for performing reparametrization in the sense of the equidistribution rule. The key idea is to regard a periodic monitor function of the arclength variable as the true curvature of an open planar curve in the canonical parametrization and to consider a linear interpolation between the curvature and that of the unit circle. This process defines a one-parameter family of planar curves, and since a simple formula is known for tangential velocity that maintains the the equidistribution rule with positive curvature, it is straightforward to find the desired correspondence between the arclength and another variable by tracking the generated ``moving" curve. With the $L^1$-type normalization of a monitor function, which periodizes derivatives of the interpolated curve, the normal and tangential velocity in the motion can be computed with spectral accuracy while the resulting parametrization remains the same. At the validation step, the whole algorithm is tested on the two prominent examples with locally high complexity, and improvements to spatial resolution are evaluated by recomputing the invariant Fourier coefficients from downsampled data. As concluding remarks, we mention three important directions for further developments. \par
Firstly, an automatic process for generating monitor functions, which are superpositions of periodized Gaussians in Section \ref{results}, is desirable in practical simulations. This problem is partly addressed in the previous work \cite{Koga2021}, where highly oscillatory signals are ``bundled" into single peaks by a Hilbert-type envelope and subsequently smoothed by the Gaussian filter. However, such a technique precludes estimations on sufficient numbers of sample points for resolving generated functions. A major advantage of defining monitor functions in terms of well-known kernels is that one can minimize the costs for discretizing interpolated curves in the reparametrization algorithm with a priori information on the fast-decaying smooth function, as found in the developments of the NUFFT algorithms \cite{BaMaKl2019, GrLe2004, KeKuPo2009}.       \par
Secondly, since a static algorithm separates mesh refinement from time-stepping, it should be possible to incorporate an adaptive stepsize control into numerical simulations of moving boundary problems. This potential improvement to the previous work \cite{Koga2021} may enable to explain, for example, experimental results such as droplet formations affected by weak viscosity \cite{DeHeHaVeRoKeEgBo2018}. Such a phenomenon is expected to be an appropriate setting for verifying improvements to the previous implementation because its local profile is visible at smaller length-scales than those of the inertia-capillary regime. As a model for this challenge, the boundary integral formulation with weak viscosity suggested by Lundgren and Mansour \cite{LuMa1988} is a promising candidate. \par
Lastly, extending the core ideas to other geometric settings is of great importance. In the context of fluid mechanics, the liquid-air interface of a two-dimensional droplet on a solid surface can be modeled by an open curve and controlling its mesh is effective in immersed boundary simulations \cite{LaTsHu2010}. Unfortunately, the Fourier-type approach in the present work is not directly applicable to the non-periodic problem, while the Chebyshev series and its coefficients should serve as analogous tools for open planar curves with acceleration by the nonuniform fast cosine transform \cite{KeKuPo2009}. On the other hand, efficient representations of surfaces diffeomorphic to the sphere, which are a natural high-dimensional analogue of periodic Jordan curves in the plane, are fundamental in boundary integral simulations of non-axisymmetric droplets, whereas practical redistribution techniques for the class are limited to spatial discretization of finite order such as simplicial or quadrilateral mesh \cite{KoHu2020, CrRuRu2015}. A major barrier to extending those methods to spectrally accurate discretization is that   parametrization of the sphere by a single smooth mapping always forms ``poles'', which is problematic not only because it is impossible to perform mesh coarsening  there but also because simulations of fluid interfaces with surface tension may break down due to stiffness \cite{AmSiTl2013}. This difficulty motivates us to first decompose the sphere into at least two parts and associate each with an appropriate region in $\mathbb{R}^2$ so that the effective approach for piecewise smooth objects \cite{KoHu2020} is naturally extended. Thus, there are many topics relevant to the present work, and we hope to address these problems in future works.


\section*{Acknowledgement}
The author would like to thank Tomoyuki Miyaji for valuable comments. We are also grateful to the two anonymous reviewers for their helpful suggestions. This work was partially supported by Kyoto University Grant for Start-Up 2021, JST ACT-X JPMJAX2106, and JSPS KAKENHI Grant No.\,21K20325.

\bibliographystyle{siamplain}
\bibliography{ref}
\end{document}